\documentclass[12pt]{amsart}
\usepackage{mathrsfs}
\usepackage{amssymb}
\usepackage{verbatim}
\usepackage{hyperref}
\usepackage{color}

\begin{document}
\newtheorem*{thmC}{Theorem C}
\newtheorem*{thmA}{Theorem A}
\newtheorem*{thmB}{Theorem B}
% define theorem environments
\newtheorem{theorem}{Theorem}    %[section]
\newtheorem{proposition}[theorem]{Proposition}
\newtheorem{conjecture}[theorem]{Conjecture}
\def\theconjecture{\unskip}
\newtheorem{corollary}[theorem]{Corollary}
\newtheorem{lemma}[theorem]{Lemma}
\newtheorem{sublemma}[theorem]{Sublemma}
\newtheorem{observation}[theorem]{Observation}
\theoremstyle{definition}
\newtheorem{definition}{Definition}
\newtheorem{notation}[definition]{Notation}
\newtheorem{remark}[definition]{Remark}
\newtheorem{question}[definition]{Question}
\newtheorem{questions}[definition]{Questions}
\newtheorem{example}[definition]{Example}
\newtheorem{problem}[definition]{Problem}
\newtheorem{exercise}[definition]{Exercise}

\numberwithin{theorem}{section}
\numberwithin{definition}{section}
\numberwithin{equation}{section}

\def\earrow{{\mathbf e}}
\def\rarrow{{\mathbf r}}
\def\uarrow{{\mathbf u}}
\def\varrow{{\mathbf V}}
\def\tpar{T_{\rm par}}
\def\apar{A_{\rm par}}

\def\reals{{\mathbb R}}
\def\torus{{\mathbb T}}
\def\heis{{\mathbb H}}
\def\integers{{\mathbb Z}}
\def\naturals{{\mathbb N}}
\def\complex{{\mathbb C}\/}
\def\distance{\operatorname{distance}\,}
\def\support{\operatorname{support}\,}
\def\dist{\operatorname{dist}\,}
\def\Span{\operatorname{span}\,}
\def\degree{\operatorname{degree}\,}
\def\kernel{\operatorname{kernel}\,}
\def\dim{\operatorname{dim}\,}
\def\codim{\operatorname{codim}}
\def\trace{\operatorname{trace\,}}
\def\Span{\operatorname{span}\,}
\def\dimension{\operatorname{dimension}\,}
\def\codimension{\operatorname{codimension}\,}
\def\nullspace{\scriptk}
\def\kernel{\operatorname{Ker}}
\def\ZZ{ {\mathbb Z} }
\def\p{\partial}
\def\rp{{ ^{-1} }}
\def\Re{\operatorname{Re\,} }
\def\Im{\operatorname{Im\,} }
\def\ov{\overline}
\def\eps{\varepsilon}
\def\lt{L^2}
\def\diver{\operatorname{div}}
\def\curl{\operatorname{curl}}
\def\etta{\eta}
\newcommand{\norm}[1]{ \|  #1 \|}
\def\expect{\mathbb E}
\def\bull{$\bullet$\ }

\def\xone{x_1}
\def\xtwo{x_2}
\def\xq{x_2+x_1^2}
\newcommand{\abr}[1]{ \langle  #1 \rangle}

\newcommand{\Norm}[1]{ \left\|  #1 \right\| }
\newcommand{\set}[1]{ \left\{ #1 \right\} }
\def\one{\mathbf 1}
\def\whole{\mathbf V}
\newcommand{\modulo}[2]{[#1]_{#2}}

\def\scriptf{{\mathcal F}}
\def\scriptg{{\mathcal G}}
\def\scriptm{{\mathcal M}}
\def\scriptb{{\mathcal B}}
\def\scriptc{{\mathcal C}}
\def\scriptt{{\mathcal T}}
\def\scripti{{\mathcal I}}
\def\scripte{{\mathcal E}}
\def\scriptv{{\mathcal V}}
\def\scriptw{{\mathcal W}}
\def\scriptu{{\mathcal U}}
\def\scriptS{{\mathcal S}}
\def\scripta{{\mathcal A}}
\def\scriptr{{\mathcal R}}
\def\scripto{{\mathcal O}}
\def\scripth{{\mathcal H}}
\def\scriptd{{\mathcal D}}
\def\scriptl{{\mathcal L}}
\def\scriptn{{\mathcal N}}
\def\scriptp{{\mathcal P}}
\def\scriptk{{\mathcal K}}
\def\frakv{{\mathfrak V}}

\begin{comment}
\def\scriptx{{\mathcal X}}
\def\scriptj{{\mathcal J}}
\def\scriptr{{\mathcal R}}
\def\scriptS{{\mathcal S}}
\def\scripta{{\mathcal A}}
\def\scriptk{{\mathcal K}}
\def\scriptp{{\mathcal P}}
\def\frakg{{\mathfrak g}}
\def\frakG{{\mathfrak G}}
\def\boldn{\mathbf N}
\end{comment}

\author{Jian Tan}
\address{Jian Tan
\\
College of Science
\\
Nanjing University of Posts and Telecommunications
\\Nanjing 210023
\\
People's Republic of China
}
\email{tanjian89@126.com}


%\thanks{The authors were supported partly by NSFC
%(No. 11471041), the Fundamental Research Funds for the Central Universities (NO. 2014kJJCA10) and NCET-13-0065. \\ \indent Corresponding author: Qingying Xue\indent Email: qyxue@bnu.edu.cn}

\title
[Multilinear fractional integrals and their commutators]
{Multilinear fractional type operators\\ and their commutators on
Hardy spaces with
variable exponents}
\maketitle

\begin{abstract}
In this article, we show that multilinear
fractional type operators are bounded
from product Hardy spaces with variable exponents
into Lebesgue spaces with variable exponents via
the atomic decomposition theory. We also study continuity properties of
commutators of multilinear fractional type operators on product of certain
Hardy spaces with variable exponents.
\end{abstract}

\section{Introduction}

The study of Hardy spaces began in the
early 1900s
in the context of Fourier series
and complex analysis in one variable.
It was not until
1960
when
the groundbreaking work in Hardy space theory in $\mathbb R^n$
came from Stein, Weiss, Coifman and C. Fefferman
in \cite{C,CWe,FS}.
The classical Hardy space can be characterized by the
Littlewood-Paley-Stein square functions, maximal functions and atomic decompositions. Especially,
atomic decomposition is a significant tool in harmonic analysis and wavelet analysis
for the study of function spaces and
the operators acting on these spaces.
Atomic decomposition was first introduced by
Coifman (\cite{C}) in one dimension in 1974 and later was extended
to higher dimensions by Latter (\cite{L}). As we all
know, atomic decompositions of Hardy spaces play an
important role in the boundedness of operators on Hardy spaces
and it is commonly sufficient to check that atoms are mapped into bounded
elements of quasi-Banach spaces.

Another stage in the progress of the theory of Hardy spaces was done by
Nakai and Sawano (\cite{NS}) and
Cruz-Uribe and Wang (\cite{CW}) recently
when they independently considered Hardy spaces
with variable exponents.
It is quiet different to obtain the boundedness of operators
on Hardy spaces with variable exponents. It is not sufficient to show
the $H^{p(\cdot)}$-boundedness
merely by checking the action
of the operators on $H^{p(\cdot)}$-atoms. In the linear theory,
the boundedness of some operators on variable Hardy spaces
have been
established in \cite{CW,Ho,NS,S,ZSY} as applications of the corresponding
atomic decompositions
theories.

In more recent years, the study of multilinear operators
on Hardy space theory has received increasing attention by many
authors, see for example \cite{GK,HM,HL}.
While the multilinear
operators worked well on the product of Hardy spaces,
it is surprising that these similar results in the setting of
variable exponents were unknown for a long time.
The boundedness of some multilinear operators
on products of classical Hardy spaces
was investigated by Grafakos and Kalton (\cite{GK})
and Li et al. (\cite{LXY}). In \cite{TLZ},
Tan et al. studied some multilinear
operators are bounded on variable Lebesgue
spaces $L^{p(\cdot)}$.
 However, there are some
subtle difficulties in proving the boundedness results when we
deal with the $H^{p(\cdot)}$-norm.
The first goal of this article
is to show that
multilinear fractional type operators
are bounded from product of Hardy spaces with variable exponents
into Lebesgue or Hardy spaces with variable exponents via
atomic decompositions theory.
We also remark that some boundedness of many types of
multilinear operators on some variable Hardy spaces have established in \cite{CMN,Tan,Tan1,TZ1}. After we were completing this paper
we learned that some of similar results
had been also established independently by Cruz-uribe et al. \cite{CMN1} independently
though our approaches are very different.
Besides, we also obtained the boundedness of commutators of multilinear fractional type operators on variable Hardy spaces, which is also interesting and useful.

On the other hand, we study the boundedness of
commutators of
multilinear fractional type operators.
In 1976, Coifman et al. (\cite{CRW})
studied the $L^p$ boundedness
of linear commutators generated by the
Calder\'{o}n-Zygmund singular integral
operator and $b\in \mbox{BMO}$.
In 1982,
Chanillo (\cite{Ch}) consider the boundedness
of commutators of fractional integral operators
on classical Lebesgue spaces.
Similar to
the property of a linear Calder\'on-Zygmund operator, a
linear fractional type operator $I_\alpha$ associated with
a BMO function $b$ fails to
satisfy the continuity from the Hardy space $H^p$ into $L^p$ for
$p\le1$. In 2002,
Ding et al. (\cite{DLZ}) proved that
$[b, I_\alpha]$ is continuous from an atomic Hardy space $H^p_b$
into $L^p$, where $H^p_b$ is a
subspace of the Hardy space $H^p$ for $n/(n+1)<p\le1$.
In addition, the boundedness of the commutators of
multilinear operators has also been studied already
in \cite{CX,M,PT}. Then, Li and Xue (\cite{LX}) consider
continuity properties for commutators of
multilinear type operators on product of
certain Hardy spaces.
It is natural to ask whether such results
are also hold in variable exponents setting.
The answer is affirmative.
The second purpose in this article is to
study the commutators of multilinear fractional type operators
on product of certain
Hardy spaces with variable exponents.
To do so, we will introduce a new atomic space with variable exponents,
$H^{p(\cdot)}_b$,
which is a subspace of the Hardy space with variable exponents
$H^{p(\cdot)}$ and obtain the endpoint
$(H_b^{p_1(\cdot)}\times\cdots\times H_b^{p_m(\cdot)}
, L^{q(\cdot)})$ boundedness for multilinear
fractional type operators.

First we recall the definition
of Lebesgue spaces with variable exponent.
Note that the variable exponent spaces,
such as the variable Lebesgue spaces and the variable
Sobolev spaces, were studied by a substantial number
of researchers (see, for instance, \cite{CFMP,KR}).
For any Lebesgue measurable function $p(\cdot):
\mathbb R^n\rightarrow (0,\infty]$ and for any
measurable subset $E\subset \mathbb{R}^n$, we denote
$p^-(E)= \inf_{x\in E}p(x)$ and $p^+(E)= \sup_{x\in E}p(x).$
Especially, we denote $p^-=p^{-}(\mathbb{R}^n)$ and $p^+=p^{+}(\mathbb{R}^n)$.
Let $p(\cdot)$: $\mathbb{R}^n\rightarrow(0,\infty)$ be a measurable
function with $0<p^-\leq p^+ <\infty$ and $\mathcal{P}^0$
be the set of all these $p(\cdot)$.
Let $\mathcal{P}$ denote the set of all measurable functions
$p(\cdot):\mathbb{R}^n \rightarrow[1,\infty) $ such that
$1<p^-\leq p^+ <\infty.$

\begin{definition}\label{s1de1}\quad
Let $p(\cdot):\mathbb R^n\rightarrow (0,\infty]$
be a Lebesgue measurable function.
The variable Lebesgue space $L^{p(\cdot)}$ consisits of all
Lebesgue measurable functions $f$, for which the quantity
$\int_{\mathbb{R}^n}|\varepsilon f(x)|^{p(x)}dx$ is finite for some
$\varepsilon>0$ and
$$\|f\|_{L^{p(\cdot)}}=\inf{\left\{\lambda>0: \int_{\mathbb{R}^n}\left(\frac{|f(x)|}{\lambda}\right)^{p(x)}dx\leq 1 \right\}}.$$
\end{definition}

The variable Lebesgue spaces
were first established by Orlicz \cite{O}
in 1931. Two decades later,
Nakano \cite{N1} first systematically
studied modular function spaces which include the variable
Lebesgue spaces as specific examples.
However, the modern development started with the paper \cite{KR}
of Kov\'{a}\v{c}ik and R\'{a}kosn\'{\i}k in 1991.
As a special case of the theory of Nakano and Luxemberg, we see that $L^{p(\cdot)}$
is a quasi-normed space. Especially, when $p^-\geq1$, $L^{p(\cdot)}$ is a Banach space.

We also recall the following class of
exponent function, which can be found in \cite{D}.
Let $\mathcal{B}$ be the set of $p(\cdot)\in \mathcal{P}$ such that the
Hardy-littlewood maximal operator $M$ is bounded on  $L^{p(\cdot)}$.
An important subset of $\mathcal{B}$ is $LH$
condition.

In the study of variable exponent function spaces it is common
to assume that the exponent function $p(\cdot)$ satisfies $LH$
condition.
We say that $p(\cdot)\in LH$, if $p(\cdot)$ satisfies

 $$|p(x)-p(y)|\leq \frac{C}{-\log(|x-y|)} ,\quad |x-y| \leq 1/2$$
and
 $$|p(x)-p(y)|\leq \frac{C}{\log|x|+e} ,\quad |y|\geq |x|.$$

It is well known
that $p(\cdot)\in \mathcal{B}$ if $p(\cdot)\in \mathcal{P}\cap LH.$
Moreover, example shows that the above $LH$ conditions are necessary in certain sense, see Pick and R\.{u}\u{z}i\u{c}ka (\cite{PR}) for more details. Next we also
recall the definition of variable Hardy spaces ${H}^{p(\cdot)}$ as follows.

\begin{definition}\label{s1de2}~(\cite{CW,NS})\quad
Let $f\in \mathcal{S'}$,
$\psi\in \mathcal S$, $p(\cdot)\in {\mathcal{P}^0}$
and $\psi_t(x)=t^{-n}\psi(t^{-1}x)$, $x\in \mathbb{R}^n$.
Denote by $\mathcal{M}$ the grand maximal operator given by
$\mathcal{M}f(x)= \sup\{|\psi_t\ast f(x)|: t>0,\psi \in \mathcal{F}_N\}$ for any fixed large integer $N$,
where $\mathcal{F}_N=\{\varphi \in \mathcal{S}:\int\varphi(x)dx=1,\sum_{|\alpha|\leq N}\sup(1+|x|)^N|\partial ^\alpha \varphi(x)|\leq 1\}$.
The variable Hardy space ${H}^{p(\cdot)}$ is the set of all $f\in \mathcal{S}^\prime$, for which the quantity
$$\|f\|_{{H}^{p(\cdot)}}=\|\mathcal{M}f\|_{{L}^{p(\cdot)}}<\infty.$$
\end{definition}

Throughout this paper, $C$ or $c$ will denote a positive constant that may vary at each occurrence
but is independent to the essential variables, and $A\sim B$ means that there are constants
$C_1>0$ and $C_2>0$ independent of the essential variables such that $C_1B\leq A\leq C_2B$.
Given a measurable set $S\subset \mathbb{R}^n$, $|S|$ denotes the Lebesgue measure and $\chi_S$
means the characteristic function. For a cube $Q$, let $Q^\ast$ denote with the same center
and $2\sqrt{n}$ its side length, i.e. $l(Q^\ast)=100\sqrt{n}l(Q)$.
The symbols $\mathcal S$ and $\mathcal S'$ denote the class of
Schwartz functions and tempered functions, respectively.
As usual, for a function $\psi$ on $\mathbb R^n$
and $\psi_t(x)=t^{-n}\psi(t^{-1}x)$.
We also use the notations $j\wedge j'=\min\{j,j'\}$ and $j\vee j'=\max\{j,j'\}$.
Moreover, denote by $L^q_{comp}$ the set of all $L^q$-functions with compact support.
For $L=0,1,2,\ldots,$ $\mathcal P_L$ denotes the set of all polynomials with degree
less than or equal to $L$ and $\mathcal P_{-1}\equiv\{0\}$.
The spaces $\mathcal P_L^{\perp}$ is the set of all integrable functions $f$
satisfying $\int_{\mathbb R^n}(1+|x|)^L|A(x)|dx<\infty$ and $\int_{\mathbb R^n}
x^\alpha A(x)dx=0$ for all multiindices $\alpha$ such that $|\alpha|\le L$.
By convention, $\mathcal P_{-1}^{\perp}$ is the set of all measurable functions.
For $L=-1,0,1,\cdots$, define $L^{q,L}_{comp}\equiv L^q_{comp}\cap \mathcal P_L^\perp$.

In what follows, we recall the new atoms for Hardy spaces with variable exponents
${H}^{p(\cdot)}$,
which is introduced in \cite{S}. Define
$$
p_-=p^-\wedge1,\quad d_{p(\cdot)}\equiv [n/{p_-}-n]\vee -1
$$
for $p\in(0,\infty).$
Let $p(\cdot): \mathbb{R}^n\rightarrow (0,\infty)$, $0< p^{-}\leq p^+\leq \infty$.
Fix an integer $d\geq d_{p(\cdot)}$ and $1<q\le\infty$.
A function $a$ on $\mathbb{R}^n$ is called a $(p(\cdot),q)$-atom, if there exists a cube $Q$
such that
${\rm supp}\,a\subset Q$;
$\|a\|_{L^q}\leq {|Q|^{1/q}}$;
$\int_{\mathbb{R}^n} a(x)x^\alpha dx=0\; {\rm for}\;|\alpha| \leq d$.
Especially, the first two conditions can be replaced by
$|a|\leq \chi_Q$ when $q=\infty$.

The atomic decomposition of Hardy spaces with variable exponents
was first established independently in \cite{CW,NS}.
Moreover, Yang et al.\cite{YYZ,YZN,ZYL} established some equivalent characterizations of Hardy spaces with variable exponents.
Recently, the
author revisited the
atomic decomposition for $H^{p(\cdot)}$ via the Littlewood-Paley-Stein
theory in \cite{T17AFA}. In this paper we will use the
following decomposition results obtained by Sawano (\cite{S}), which
extends and sharp the ones of above papers.

\begin{theorem}\label{s1th01}\cite{S}\quad Let $p(\cdot)\in LH\cap \mathcal P^0$
and $q>(p^+\vee 1)$.
Suppose that $d\geq d_{p(\cdot)}$ and $s\in(0,\infty)$.
If $f\in H^{p(\cdot)}$, there exists
sequences of $(p(\cdot),\infty)-$atoms $\{a_j\}$ and scalars $\{\lambda_j\}$
{\noindent}such that $f=\sum_{j=1}^\infty\lambda_ja_j$
in $H^{p(\cdot)}\cap L^q$
and that
\begin{equation*}
 \left\|\left\{\sum_{j=1}^\infty\left({\lambda_j
 \chi_{Q_j}}\right)^{s}
\right\}^{\frac{1}{s}}\right\|_{L^{p(\cdot)}}
  \leq C\|f\|_{{H}^{p(\cdot)}}.
\end{equation*}
\end{theorem}

The multilinear fractional type operators are natural
generalization of linear ones.
Their earliest version was originated on the work
of Grafakos (\cite{G}) in 1992, in which he studied the multilinear fractional integral
defined by
\begin{align*}
\bar{I}_\alpha(\vec{f})(x)=\int_{\mathbb R^n}
\frac{1}{|y|^{n-\alpha}}\prod_{i=1}^mf_i(x-\theta_iy)dy,
\end{align*}
where $\theta_i(i=1,\cdots,m)$ are fixed distinct
and nonzero real numbers and $0<\alpha<n.$
Later on, in 1998, Kenig and Stein \cite{KS} established
the boundedness of another
type of multilinear fractional integral
${I}_{\alpha,A}$
on product of Lebesgue spaces.
${I}_{\alpha,A}(\vec{f})$ is defined by
\begin{align*}
{I}_{\alpha,A}(\vec{f})(x)=\int_{(\mathbb R^n)^m}
\frac{1}{|(y_1,\ldots,y_m)|^{mn-\alpha}}\prod_{i=1}^mf_i(
\ell_i(y_1,\ldots,y_m,x))dy_i,
\end{align*}
where $\ell_i$ is a linear combination of $y_j$＊s and $x$
depending on the matrix $A$.
In \cite{LL}, Lin and Lu obtained $I_{\alpha,A}$ is bounded
from product of Hardy spaces to Lebegue spaces when
$\ell_i(y_1, \ldots, y_m, x)=x-y_i$. We denote this multilinear fractional
type integral operators by $I_\alpha$, namely,
\begin{align*}
{I}_{\alpha}(\vec{f})(x)=\int_{(\mathbb R^n)^m}
\frac{1}{|(y_1,\ldots,y_m)|^{mn-\alpha}}\prod_{i=1}^mf_i
(x-y_i)dy_i.
\end{align*}
For convenience, we also denote
$K_\alpha(y_1,\ldots,y_m)=\frac{1}{|(y_1,\ldots,y_m)|^{mn-\alpha}}$.

For any $1\le j\le m$, we can define the commutator of multilinear operator by
$$[b,T]_j(\vec{f})(x):=bT(\vec{f})(x)-T(f_1,\ldots,bf_j,\ldots,f_m)(x),$$
where $b$ is a locally integral function and $T$ is a multilinear
operator.

Then $[b,I_\alpha]_j(\vec{f})$ is defined by
$$[b,I_\alpha]_j(\vec{f})(x)
=\int_{(\mathbb{R}^n)^m}\frac{(b(x)-b(y_j))\prod_{i=1}^{m}f_i(y_i)}
{(\sum_{i=1}^{m}|x-y_i|)^{nm-\alpha}}\prod_{i=1}^mdy_i.$$

In Section 2, we will show that ${I}_{\alpha}$
is bounded from product of Hardy spaces with variable exponents
to Lebegue or Hardy spaces with exponents. Then we will introduce
the new atomic Hardy spaces with variable exponents
$H^{p(\cdot)}_b$ in section 3. Moreover,
we also consider continuity properties for commutators of
multilinear type operators on product of
the atomic Hardy spaces with variable exponents $H^{p(\cdot)}_b$.

\section{Multilinear fractional type operators
on product of Hardy spaces with variable exponents}

In this section, we will discuss the boundedness of
multilinear fractional type operators on product of Hardy spaces
with variable exponents. The results are new even
of the classsical constant $\prod_{j=1}^m H^{p_j}\rightarrow H^{q}$
boundedness for multilinear fractional type operators.
First we introduce some necessary
notations and requisite lemmas.
The following generalized H\"{o}lder inequality on variable Lebesgue spaces
can be found in in \cite{CF} or \cite{TLZ}.

\begin{lemma}\label{s2l1}
\quad Given exponent function $p_i(\cdot)\in \mathcal{P}^0,$ define
$p(\cdot)\in \mathcal{P}^0$ by
$$\frac{1}{p(x)}=\sum_{i=1}^m\frac{1}{p_i(x)},$$
where $i=1,\ldots,m.$
Then for all $f_i\in L^{p_i(\cdot)}$ and
$f_1\cdots f_m\in L^{p(\cdot)}$ and
$$\|\prod_{i=1}^mf_i\|_{p(\cdot)}\leq C\prod_{i=1}^m\|f_i\|_{p_i(\cdot)}.$$
\end{lemma}

\begin{lemma}(\cite{TLZ})\label{s2l2}\quad Let $m\in\mathbb N$,
$$\frac{1}{s(x)}=\sum_{i=1}^m\frac{1}{r_i(x)}-\frac{\alpha}{n}, x\in\mathbb{R}^n,$$
with $0<\alpha<mn$, $1< r_i\le\infty$. Then
\begin{equation*}
  \|I_\alpha(\vec{f})\|_{s(\cdot)}\leq C\prod_{i=1}^m\|f_i\|_{r_i(\cdot)}.
\end{equation*}
\end{lemma}

We also need the following boundedness of the
vector-valued fractional maximal
operators on variable Lebesgue spaces whose proof can be found in
\cite{CF}.
Let $0\leq\alpha<n$, we define $$M_\alpha f(x)=\sup_{Q\ni x}\frac{1}{|Q|^{1-\alpha/n}}\int_Q|f(y)|dy.$$

\begin{lemma}\label{s2l3}\quad Let $0\leq\alpha<n$, $p(\cdot),q(\cdot)\in \mathcal{B}$
be such that $p^+<\frac{n}{\alpha}$ and
$$\frac{1}{p(x)}-\frac{1}{q(x)}=\frac{\alpha}{n},\quad x\in\mathbb{R}^n.$$
If $q(\cdot)(n-\alpha)/n\in \mathcal{B}$, then
for any $q>1$, $f=\{f_i\}_{i\in \mathbb{Z}}$, $f_i\in L_{loc}$, $i\in \mathbb{Z}$
\begin{equation*}
  \|\|\mathbb{M_\alpha}(f)\|_{l^q}\|_{L^{q(\cdot)}}\leq C\|\|f||_{l^q}\|_{L^{p(\cdot)}},
\end{equation*}
where $\mathbb{M_\alpha}(f)=\{M_\alpha(f_i)\}_{i\in\mathbb{Z}}$.
\end{lemma}

\begin{lemma}\cite{KR}\label{s2l4} \quad Let $p(\cdot)\in \mathcal{P}$,
$f\in L^{p(\cdot)}$ and $g\in L^{p'(\cdot)}$, then
$fg$ is integrable on $\mathbb{R}^n$ and
$$\int_{\mathbb{R}^n}|f(x)g(x)|dx\leq r_p\|f\|_{L^{p(\cdot)}}
\|g\|_{L^{p'(\cdot)}},$$
where $r_p= 1+1/p^--1/p^+.$
Moreover, for all $g\in L^{p'(\cdot)}$ such that $\|g\|_{L^{p'(\cdot)}}\le1$ we get that
$$
\|f\|_{L^{p(\cdot)}}\le \sup_{g}|\int_{\mathbb R^n}f(x)g(x)dx|
\le r_p\|f\|_{L^{p(\cdot)}}.
$$
\end{lemma}

Our first theorem is the following.
\begin{theorem}\label{s1th1}\quad
Let $0<\alpha<n$.
Suppose that $p_1(\cdot),\ldots,p_m(\cdot)
\in LH\cap \mathcal P^0$
and $q(\cdot)\in \mathcal P^0$ be Lebesgue measure functions satisfying
\begin{align}\label{s1c2}
\frac{1}{p_1(x)}+\cdots+\frac{1}{p_m(x)}-\frac{\alpha}{n}=\frac{1}{q(x)},
\quad\quad x\in\mathbb R^n.
\end{align}
Then $I_\alpha$ can be extended to a bounded operator from
$\prod_{j=1}^m H^{p_j(\cdot)}$ into $L^{q(\cdot)}$.
\end{theorem}

\noindent\textit{\bf Proof}\quad
Observe that $p_j(\cdot)\in LH\cap \mathcal P^0$
and choose that $\bar q>(p_j^+\vee 1)$, $j=1,\cdots,m$.
By Theorem \ref{s1th01}, for each
$f_j\in{H}^{p_j(\cdot)}\cap{L}^{\bar q}$, $j=1,\ldots,m$,
$f_j$ admits an atomic decomposition:
Suppose that $d_j\geq d_{p_j(\cdot)}$ and $s\in(0,\infty)$.
If $f\in H^{p_j(\cdot)}$, there exists
sequences of $(p_j(\cdot),\infty)-$atoms $\{a_{j,k_j}\}$
and scalars $\{\lambda_{j,k_j}\}$
{\noindent}such that $f_j=\sum_{k_j=1}^\infty\lambda_{j,k_j}a_{j,k_j}$
in $H^{p_j(\cdot)}\cap L^{\bar q}$
and that
\begin{equation*}
 \left\|\left\{\sum_{k_j=1}^\infty\left({\lambda_{j,k_j}
 \chi_{Q_{j,k_j}}}\right)^{s}
\right\}^{\frac{1}{s}}\right\|_{L^{p_j(\cdot)}}
  \leq C\|f_j\|_{{H}^{p_j(\cdot)}}.
\end{equation*}

For the decomposition of $f_j$, $j=1,\ldots,m$,
we write
\begin{align*}\label{}
I_\alpha(\vec{f})(x)=\sum_{k_1}\cdots\sum_{k_m}\lambda_{1,k_1}
\cdots\lambda_{m,k_m}I_\alpha(a_{1,k_1},\ldots,a_{m,k_m})(x)
\end{align*}
in the sense of distributions.

Fixed $k_1,\cdots,k_m$, there are two cases for $x\in\mathbb R^n$.\\
Case 1: $x\in Q^\ast_{1,k_1}\cap\cdots\cap Q^\ast_{m,k_m}$.\\
Case 2: $x\in Q^{\ast,c}_{1,k_1}\cup\cdots\cup Q^{\ast,c}_{m,k_m}$.\\

Then we have
\begin{align*}\label{}
|I_\alpha(\vec{f})(x)|&\le\sum_{k_1}\cdots\sum_{k_m}|\lambda_{1,k_1}|
\cdots|\lambda_{m,k_m}||I_\alpha(a_{1,k_1},\ldots,a_{m,k_m})(x)|
\chi_{Q^\ast_{1,k_1}\cap\cdots\cap Q^\ast_{m,k_m}}(x)\\
+&\sum_{k_1}\cdots\sum_{k_m}|\lambda_{1,k_1}|
\cdots|\lambda_{m,k_m}||I_\alpha(a_{1,k_1},\ldots,a_{m,k_m})(x)|
\chi_{Q^{\ast,c}_{1,k_1}\cup\cdots\cup Q^{\ast,c}_{m,k_m}}(x)\\
&=I_1(x)+I_2(x).
\end{align*}

First, we consider the estimate of $I_1(x)$.
We will show that
\begin{align}\label{I1}
\begin{split}
\|I_1\|_{L^{q(\cdot)}}
\leq C\prod_{j=1}^m\|f_j\|_{{H}^{p_j(\cdot)}}.
\end{split}
\end{align}
For fixed $k_1,\ldots,k_m$,
assume that $Q^\ast_{1,k_1}\cap\cdots\cap Q^\ast_{m,k_m}\neq 0$, otherwise
there is nothing to prove. Without loss of generality,
suppose that $Q_{1,k_1}$ has the smallest size among all these cubes.
We can pick a cube $G_{k_1,\cdots,k_m}$ such that
$$Q^\ast_{1,k_1}\cap\cdots\cap Q^\ast_{m,k_m}\subset
G_{k_1,\ldots,k_m}\subset G^\ast_{k_1,\ldots,k_m}\subset
Q^{\ast\ast}_{1,k_1}\cap\cdots\cap Q^{\ast\ast}_{m,k_m}$$
and $|G_{k_1,\ldots,k_m}|\ge C|Q_{1,k_1}|$.

Denote $H(x):=|I_\alpha(a_{1,k_1},\ldots,a_{m,k_m})(x)|
\chi_{Q_{1,k_1}^\ast\cap\cdots\cap Q_{m,k_m}^\ast}(x)$.
Obviously, $$\mbox{supp} H(x)\subset
{Q_{1,k_1}^\ast\cap\cdots\cap Q_{m,k_m}^\ast}\subset
G_{k_1,\cdots,k_m}.$$

By Lemma \ref{s2l2},
since $I_\alpha$ maps $L^{r_1}\times L^\infty\times\cdots
\times L^\infty$ into $L^{s}$ ($s>1$ and
$\frac{1}{s}=\frac{1}{r_1}-\frac{\alpha}{n}$), we
get that
\begin{align}\label{s2i3}
\begin{split}
\|H\|_{L^{s}}&\le
\|I_\alpha(a_{1,k_1},\ldots,a_{m,k_m})\|_{L^{s}}\\
&\le
C\|a_{1,k_1}\|_{L^{r_1}}\|a_{2,k_2}\|_{L^\infty}
\cdots\|a_{m,k_m}\|_{L^\infty}\\
&\le C|Q_{1,k_1}|^{\frac{1}{r_1}}
\|a_{2,k_2}\|_{L^\infty}\|a_{m,k_m}\|_{L^\infty}\\
&\le C
|Q_{1,k_1}|^{\frac{1}{r_1}}\\
&\le C
|G_{k_1,\ldots,k_m}|^{\frac{1}{r_1}}.
\end{split}
\end{align}

For any $g\in L^{(q(\cdot)/{q^-})'}$ with
$\|g\|_{L^{(q(\cdot)/{q^-}})'}\le 1$,
by H\"older inequality and (\ref{s2i3})
we find that
\begin{align*}
&\left|\int_{\mathbb R^n}H(x)^{q^-}g(x)dx\right|
\le \|H^{q^-}\|_{L^{s/q^-}}\|g\|_{L^{(s/q)'}}\\
\le& C
|G_{k_1,\ldots,k_m}|^{\frac{q^-}{r_1}}
 \left(\int_{G_{k_1,\ldots,k_m}}|g(x)|^{(s/q^-)'}dx\right)^{\frac{1}{(s/q^-)'}},
\end{align*}
where $(s/q^-)'$ is the conjugate of $s/q^-$. Thus,
\begin{align*}
&\left|\int_{\mathbb R^n}H(x)^{q^-}g(x)dx\right|\\
\le& C
|G_{k_1,\ldots,k_m}|^{1+\frac{\alpha q^-}{n}}
 \left(\frac{1}{|G_{k_1,\ldots,k_m}|}\int_{G_{k_1,\ldots,k_m}}
 |g(x)|^{(s/q^-)'}dx\right)^{\frac{1}{(s/q^-)'}}\\
\le& C
|G_{k_1,\ldots,k_m}|^{\frac{\alpha q^-}{n}}
|G_{k_1,\ldots,k_m}|
 \inf_{x\in R_{k_1,k_2}}
 M(|g|^{(s/q^-)'})(x)^{\frac{1}{(s/q^-)'}}\\
\le& C
|G_{k_1,\ldots,k_m}|^{\frac{\alpha q^-}{n}}
\int_{G_{k_1,\ldots,k_m}}(M(|g|^{(s/q^-)'}))^{\frac{1}{(s/q^-)'}}(x)dx.
\end{align*}

When $0<q^-\le1$, using Lemma \ref{s2l4} we obtain that
\begin{align*}
\Bigg|\int_{\mathbb R^n}&\left(\sum_{k_1,\cdots,k_m}|
\prod_{j=1}^m\lambda_{j,k_j}|^{q^-}
|I_\alpha(a_{1,k_1},\ldots,a_{m,k_m})(x)|^{q^-}
\chi_{Q_{1,k_1}^\ast\cap\cdots\cap Q_{m,k_m}^\ast}(x)\right)g(x)dx
\Bigg|\\
\le&
\sum_{k_1}\cdots\sum_{k_m}|\prod_{j=1}^m\lambda_{j,k_j}|^{q^-}
\left|\int_{G_{k_1,\ldots,k_m}}
H(x)^{q^-}g(x)dx\right|\\
\le& C
\sum_{k_1}\cdots\sum_{k_m}|
\prod_{j=1}^m\lambda_{j,k_j}|^{q^-}
|G_{k_1,\ldots,k_m}|^{\frac{\alpha q^-}{n}}\\
&\times\int_{G_{k_1,\ldots,k_m}}(M(|g|^{(s/q^-)'}))^{\frac{1}{(s/q^-)'}}(x)dx\\
=&C
\int_{\mathbb R^n}\sum_{k_1}\cdots\sum_{k_m}
\prod_{j=1}^m
|\lambda_{j,k_j}|^{q^-}\chi_{G_{k_1,\ldots,k_m}}
|G_{k_1,\ldots,k_m}|^{\frac{\alpha q^-}{n}}\\
&\times(M(|g|^{(s/q^-)'}))^{\frac{1}{(s/q^-)'}}(x)dx\\
\le& C
\left\|\sum_{k_1}\cdots\sum_{k_m}
\prod_{j=1}^m
|\lambda_{j,k_j}|^{q^-}
|G_{k_1,\ldots,k_m}|^{\frac{\alpha q^-}{n}}
\chi_{G_{k_1,\ldots,k_m}}
\right\|_{L^{q(\cdot)/{q^-}}}\\
&\times\|(M(|g|^{(s/q^-)'}))^{\frac{1}{(s/q^-)'}}\|_{L^{(q(\cdot)/{q^-}})'}\\
\le& C
\left\|\sum_{k_1}\cdots\sum_{k_m}
\prod_{j=1}^m
|\lambda_{j,k_j}|^{q^-}
|G_{k_1,\ldots,k_m}|^{\frac{\alpha q^-}{n}}
\chi_{G_{k_1,\ldots,k_m}}
\right\|_{L^{q(\cdot)/{q^-}}}\\
&\times
\|(M(|g|^{(s/q^-)'}))\|^{\frac{1}{(s/q^-)'}}_{L^{(q(\cdot)/{q^-})'/{(s/q^-)'}}}.
\end{align*}

Choose $s$ large enough such that $(q(\cdot)/{q^-})'/{(s/q^-)'}>1$.
Then by Hardy-Littlewood operator $M$ is bounded on
$L^{q'(\cdot)/{(s/q^-)'q^-}}$ and $\|g\|_{L^{(q(\cdot)/{q^-})'}}\le 1$,
we know that
$$
\|(M(|g|^{(s/q^-)'}))\|^{\frac{1}{(s/q^-)'}}_{L^{(q(\cdot)/{q^-})'/{(s/q^-)'}}}
\le C.
$$

Applying Lemma \ref{s2l4} again,
we get that
\begin{align*}\label{}
\begin{split}
&\Bigg\|\sum_{k_1}\cdots\sum_{k_m}|
\lambda_{1,k_1}|^{q^-}\cdots|\lambda_{m,k_m}|^{q^-}
|I_\alpha(a_{1,k_1},\ldots,a_{m,k_m})(x)|^{q^-}
\chi_{Q_{1,k_1}^\ast\cap\cdots\cap Q_{m,k_m}^\ast}\Bigg\|_{L^{q(\cdot)/{q^-}}}\\
&\le C
\left\|\sum_{k_1}\cdots\sum_{k_m}
\prod_{j=1}^m
|\lambda_{j,k_j}|^{q^-}
|G_{k_1,\ldots,k_m}|^{\frac{\alpha q^-}{n}}
\chi_{G_{k_1,\ldots,k_m}}
\right\|_{L^{q(\cdot)/{q^-}}}\\
&\le
C\left\|\prod_{j=1}^m
\sum_{k_j}
|\lambda_{j,k_j}|^{q^-}
|Q_{j,k_j}|^{\frac{\alpha q^-}{mn}}
\chi_{Q^{\ast\ast}_{j,k_j}}
\right\|_{L^{q(\cdot)/{q^-}}}.
\end{split}
\end{align*}

Denote $\frac{1}{q_j(x)}=\frac{1}{p_j(x)}-\frac{\alpha}{mn}$
for any $x\in\mathbb R^n$, $j=1,\ldots,m$.
Then $q_j(\cdot)\in \mathcal P^0$ and for any $x\in\mathbb R^n$
$$
\frac{1}{q(x)}=\frac{1}{p_1(x)}+\cdots+\frac{1}{p_m(x)}-\frac{\alpha}{n}
=\frac{1}{q_1(x)}+\cdots+\frac{1}{q_m}.
$$

By Lemma \ref{s2l1}, we obtain
\begin{align}\label{I11}
\begin{split}
\left\|\prod_{j=1}^m
\sum_{k_j}
|\lambda_{j,k_j}|^{q^-}
|Q_{j,k_j}|^{\frac{\alpha q^-}{mn}}
\chi_{Q^{\ast\ast}_{j,k_j}}
\right\|_{L^{q(\cdot)/{q^-}}}^{\frac{1}{q^-}}
&\le \prod_{j=1}^m\left\|
\sum_{k_j}
|\lambda_{j,k_j}|^{q^-}
|Q_{j,k_j}|^{\frac{\alpha q^-}{mn}}
\chi_{Q^{\ast\ast}_{j,k_j}}
\right\|_{L^{q_j(\cdot)/{q^-}}}^{\frac{1}{q^-}}
\end{split}
\end{align}

Furthermore, it is easy to verify that
\begin{equation*}
  |Q_{j,k_j}|^{\frac{\alpha}{mn}}\chi_{Q^{\ast\ast}_{j,k_j}}(x)\leq CM_{{\alpha}{q^-}/2m}(\chi_{Q^{\ast\ast}_{j,k_j}})^{\frac{2}{q^-}}(x).
\end{equation*}

Applying Lemma 2.6, then we get that
\begin{align}\label{I12}
\begin{split}
&\left\|
\sum_{k_j}
|\lambda_{j,k_j}|^{q^-}
|Q_{j,k_j}|^{\frac{\alpha q^-}{mn}}
\chi_{Q^{\ast\ast}_{j,k_j}}
\right\|_{L^{q_j(\cdot)/{q^-}}}^{\frac{1}{q^-}}
=
\left\|
\left(\sum_{k_j}
\left(|\lambda_{j,k_j}|
|Q_{j,k_j}|^{\frac{\alpha }{mn}}
\chi_{Q^{\ast\ast}_{j,k_j}}\right)^{q^-}\right)^{\frac{1}{q^-}}
\right\|_{L^{q_j(\cdot)}}\\
&\le C\left\|
\left(\sum_{k_j}
\left(|\lambda_{j,k_j}|
M_{{\alpha}{q^-}/2m}(\chi_{Q^{\ast\ast}_{j,k_j}})^{\frac{2}{q^-}}
\right)^{q^-}\right)^{\frac{1}{q^-}}
\right\|_{L^{q_j(\cdot)}}
=
C\left\|
\left(\sum_{k_j}
|\lambda_{j,k_j}|^{q^-}
M_{{\alpha}{q^-}/2m}(\chi_{Q^{\ast\ast}_{j,k_j}})^2
\right)^{\frac{1}{q^-}}
\right\|_{L^{q_j(\cdot)}}\\
&\le
C\left\|
\left(\sum_{k_j}
|\lambda_{j,k_j}|^{q^-}
M_{{\alpha}{q^-}/2m}(\chi_{Q^{\ast\ast}_{j,k_j}})^2
\right)^{\frac{1}{2}}
\right\|^{\frac{2}{q^-}}_{L^{2q_j(\cdot)/{q^-}}}
\le
C\left\|
\left(\sum_{k_j}
|\lambda_{j,k_j}|^{q^-}
\chi_{Q^{\ast\ast}_{j,k_j}}
\right)^{\frac{1}{2}}
\right\|^{\frac{2}{q^-}}_{L^{2p_j(\cdot)/{q^-}}}.
\end{split}
\end{align}

Applying Fefferman-Stein vector value inequality on $L^{2p_j(\cdot)/{q^-}}$,
we get that
\begin{align}\label{I13}
\begin{split}
&\left\|
\left(\sum_{k_j}
|\lambda_{j,k_j}|^{q^-}
\chi_{Q^{\ast\ast}_{j,k_j}}
\right)^{\frac{1}{2}}
\right\|^{\frac{2}{q^-}}_{L^{2p_j(\cdot)/{q^-}}}
\le
\left\|
\left(\sum_{k_j}
|\lambda_{j,k_j}|^{q^-}
M{(\chi_{Q_{j,k_j}})^2}
\right)^{\frac{1}{2}}
\right\|^{\frac{2}{q^-}}_{L^{2p_j(\cdot)/{q^-}}}\\
&\le C\left\|
\left(\sum_{k_j}
\left(|\lambda_{j,k_j}|
\chi_{Q^{\ast\ast}_{j,k_j}}\right)^{q^-}
\right)^{\frac{1}{q^-}}
\right\|_{L^{p_j(\cdot)}}\le \|f\|_{H^{p_j(\cdot)}}
\end{split}
\end{align}

Therefore,
when $0<q^-\le1$,
by (\ref{I11}), (\ref{I12}) and (\ref{I13}) we have that
\begin{align*}\label{}
\begin{split}
&\|I_1\|_{L^{q(\cdot)}}
=\|(I_1)^{q^-}\|_{L^{q(\cdot)/{q^-}}}^{\frac{1}{q^-}}\\
&\le C\left\|\prod_{j=1}^m
\sum_{k_j}
|\lambda_{j,k_j}|^{q^-}
|Q_{j,k_j}|^{\frac{\alpha q^-}{mn}}
\chi_{Q^{\ast\ast}_{j,k_j}}
\right\|_{L^{q(\cdot)/{q^-}}}^{\frac{1}{q^-}}\\
&\le C\prod_{j=1}^m \|f\|_{H^{p_j(\cdot)}}
\end{split}
\end{align*}

When $q^->1$, repeating similar but more easier argument, we can
also get the desired result (\ref{I1}). In fact, we only need to
replace $p^-$ by $1$ in the above proof. We omit the detail.

Secondly, we consider the estimate of $I_2$.
Let A be a nonempty
subset of $\{1,\ldots,m\}$, and we denote the cardinality of
$A$ by $|A|$, then $1\le|A|\le m$. Let $A^c=\{1,\ldots,m\}\backslash
A.$ If $A=\{1,\ldots,m\}$, we define
$$(\cap_{j\in A}Q_{j,k_j}^{\ast,c})
\cap (\cap_{j\in A^c}Q_{j,k_j}^{\ast,c})=
\cap_{j\in A}Q_{j,k_j}^{\ast,c},$$
then
$$Q_{1,k_1}^{\ast,c}\cup\cdots\cup Q_{m,k_m}^{\ast,c}
=\cup_{A\subset\{1,\ldots,m\}}((\cap_{j\in A}Q_{j,k_j}^{\ast,c})
\cap (\cap_{j\in A^c}Q_{j,k_j}^{\ast})).$$

Set $E_A=(\cap_{j\in A}Q_{j,k_j}^{\ast,c})
\cap (\cap_{j\in A^c}Q_{j,k_j}^{\ast}).$
For fixed $A$, assume that $Q_{\tilde j,k_{\tilde j}}$
is the smallest cubes in the set of cubes $Q_{j,k_j},
j\in A$.
Let $z_{\tilde j,k_{\tilde j}}$ is the center of the cube
$Q_{\tilde j,k_{\tilde j}}$.

Denote $K_\alpha(x,y_1,\ldots,y_m)
=|(x-y_1,\ldots,x-y_m)|^{-mn+\alpha}$.
Notice that for all $|\beta|=d+1,\;
\beta=(\beta_1,\ldots,\beta_m)$
\begin{align}\label{I21}
\begin{split}
|\partial_{y_1}^{\beta_1}\cdots\partial_{y_m}^{\beta_m}
K_\alpha(x,y_1,\ldots,y_m)|
\le C|(x-y_1,\ldots,x-y_m)|^{-mn+\alpha-|\beta|}.
\end{split}
\end{align}

Since $a_{\tilde j,k_{\tilde j}}$
has zero vanishing moment up to order $d_j$,
using Taylor expansion we get
\begin{align*}
\begin{split}
  &I_\alpha(a_{1,k_1},\ldots,a_{m,k_m})(x) \\
  =& \int_{(\mathbb R^n)^m}K_\alpha(x,y_1,\ldots,y_m)
a_{1,k_1}(y_1)\cdots a_{m,k_m}(y_m)d\vec{y}\\
=& \int_{(\mathbb R^n)^{m-1}}\prod_{j\ne \tilde j}a_{j,k_j}(y_j)
\int_{\mathbb R^n}[K_\alpha(x,y_1,\ldots,y_m)
-P_{z_{\tilde j,k_{\tilde j}}}^d(x,y_1,\ldots,y_m)]
a_{\tilde j,k_{\tilde j}}d\vec{y}\\
=& \int_{(\mathbb R^n)^{m-1}}\prod_{j\ne \tilde j}a_{j,k_j}(y_j)
\int_{\mathbb R^n}\sum_{|\gamma|=d+1}
(\partial_{y_{\tilde j}}^\gamma
K_\alpha)(x,y_1,\ldots,\xi,\ldots,y_m)\frac{(y_{\tilde j}
-z_{\tilde j,k_{\tilde j}})^\gamma}{\gamma!}
a_{\tilde j}(y_{\tilde j})d\vec{y}
\end{split}
\end{align*}
for some $\xi$ on the line segment joining
$y_{\tilde j}$ to $z_{\tilde j,k_{\tilde j}}$,
where $P_{z_{\tilde j,k_{\tilde j}}}^d(x,y_1,\ldots,y_m)$
is Taylor polynomial of
$K_\alpha(x,y_1,\ldots,y_m)$.
Since $x\in (Q_{\tilde j,k_{\tilde j}}^\ast)^c$,
we can easily obtain that
$|x-\xi|\ge \frac{1}{2}|x-z_{\tilde j,k_{\tilde j}}|$.
Similarly,
$|x-y_j|\ge \frac{1}{2}|x-z_{j,k_{j}}|$
for $y_j\in Q_{j,k_j}$, $j\in A\backslash\{\tilde j\}$.

Applying the estimate for the kernel $K_\alpha$ satisfies
(\ref{I21})
and the size estimates for the new atoms yield

\begin{align*}
\begin{split}
& \int_{(\mathbb R^n)^{m-1}}\prod_{j\ne \tilde j}|a_{j,k_j}(y_j)|
\int_{\mathbb R^n}\sum_{|\gamma|=d+1}
|(\partial_{y_{\tilde j}}^\gamma
K_\alpha)(x,y_1,\ldots,\xi,\ldots,y_m)|\frac{|y_{\tilde j}
-z_{\tilde j,k_{\tilde j}}|^\gamma}{\gamma!}
|a_{\tilde j}(y_{\tilde j})|d\vec{y}\\
\le&
C\int_{(\mathbb R^n)^{|A|}}\prod_{j\in A}|a_{j,k_j}(y_j)|
\int_{(\mathbb R^n)^{m-|A|}}
\frac{|y_{\tilde j}-z_{\tilde j,k_{\tilde j}}|^{d+1}}
{(|x-\xi|+\sum_{j\ne \tilde j}|x-y_j|)^{mn+d+1-\alpha}}
\prod_{j\in A^c}|a_{ j,k_{ j}}(y_{\tilde j})|d\vec{y}\\
\le&
C\bigg(\prod_{j\in A}\|a_{j,k_j}\|_{L^1}\bigg)
\bigg(\prod_{j\in A^c}\|a_{ j,k_{ j}}\|_{L^\infty}\bigg)
\int_{(\mathbb R^n)^{m-|A|}}
\frac{|y_{\tilde j}-z_{\tilde j,k_{\tilde j}}|^{d+1}}
{(|x-\xi|+\sum_{j\ne \tilde j}|x-y_j|)^{mn+d+1-\alpha}}
d\vec{y}_{A^c}\\
\le&
C\bigg(\prod_{j\in A}\|a_{j,k_j}\|_{L^1}\bigg)
\bigg(\prod_{j\in A^c}\|a_{ j,k_{ j}}\|_{L^\infty}\bigg)
\frac{|y_{\tilde j}-z_{\tilde j,k_{\tilde j}}|^{d+1}}
{(\sum_{j\in A}|x-z_{j,k_j}|)^{mn+d+1-\alpha-n(m-|A|)}}\\
\le&
C\bigg(\prod_{j\in A}{|Q_{j,k_j}|}
\bigg)
\bigg(\prod_{j\in A^c}
{\|a_{ j,k_{ j}}\|_{L^\infty}}\bigg)
\frac{|Q_{\tilde j,k_{\tilde j}}|^{(d+1)/n}}
{(\sum_{j\in A}|x-z_{j,k_j}|)^{mn+d+1-\alpha-n(m-|A|)}}.
\end{split}
\end{align*}

Observe that $x\in \cap_{j\in A}Q_{j,k_j}^{\ast,c}$,
then we can found constant $C$ such that
$|x-z_{j,k_j}|\ge C(|x-z_{j,k_j}|+l{(Q_{j,k_j})})$.
On the other hand, using the fact that
$x\in \cap_{j\in A^c}Q_{j,k_j}^\ast$
yields that there exists a constant $C$ such that
$|x-z_{j,k_j}|\le Cl{(Q_{j,k_j})}$ for $j\in A^c$. Then
we have that
\begin{align*}
\frac{|Q_{j,k_j}|^{1+\frac{d+1}{n|A|}}}
{(|x-z_{j,k_j}|+l(Q_{j,k_j}))^{n+\frac{d+1}{|A|}}}
\ge C, \quad \mbox{for}\quad j\in A^c.
\end{align*}

Moreover, since $Q_{\tilde j,k_{\tilde j}}$
is the smallest cube among $\{Q_{j,l_j}\}_{j\in A}$,
we have that $$|Q_{\tilde j,k_{\tilde j}}|\le
\prod_{j\in A}|Q_{j,l_j}|^{\frac{1}{|A|}}.$$
Thus,
\begin{align}\label{I20}
\begin{split}
&|I_\alpha(a_{1,k_1},\ldots,a_{m,k_m})(x)|\\
\le&
C\bigg(\prod_{j\in A}\frac{|Q_{j,k_j}|^{1+\frac{d+1}{n|A|}}}
{(|x-z_{j,k_j}|+l(Q_{j,k_j}))^{n+\frac{d+1}{|A|}-\frac{\alpha}{|A|}}}\bigg)
\bigg(\prod_{j\in A^c}
{\|a_{ j,k_{ j}}\|_{L^\infty}}\bigg)\\
\le&
C\prod_{j\in A}\frac{|Q_{j,k_j}|^{1+\frac{d+1}{n|A|}}}
{(|x-z_{j,k_j}|+l(Q_{j,k_j}))^{n+\frac{d+1}{|A|}-\frac{\alpha}{|A|}}}
\prod_{j\in A^c}\frac{|Q_{j,k_j}|^{1+\frac{d+1}{n|A|}}}
{(|x-z_{j,k_j}|+l(Q_{j,k_j}))^{n+\frac{d+1}{|A|}}}
\end{split}
\end{align}
for all $x\in E_A$.

Then applying the generalized H\"older's inequality
in variable Lebesgue
spaces and Fefferman-Stein inequality in Lemma
\ref{s2l3} we obtain the estimate
\begin{align}\label{I22}
\begin{split}
&\|I_2\|_{L^{q(\cdot)}}
\le
C\Bigg\|\sum_{k_1}\cdots\sum_{k_m}\prod_{j=1}^m|\lambda_{j,{k_j}}|
\sum_{A\subset\{1,\ldots,m\}}\prod_{j\in A}\frac{|Q_{j,k_j}|^{1+\frac{d+1}{n|A|}}}
{(|x-z_{j,k_j}|+l(Q_{j,k_j}))^{n+\frac{d+1}{|A|}-\frac{\alpha}{|A|}}}
\\
&\times
\prod_{j\in A^c}\frac{|Q_{j,k_j}|^{1+\frac{d+1}{n|A|}}}
{(|x-z_{j,k_j}|+l(Q_{j,k_j}))^{n+\frac{d+1}{|A|}}}
\chi_{E_A}\Bigg\|_{L^{q(\cdot)}}\\
\le&
C\sum_{A\subset\{1,\ldots,m\}}\Bigg\|
\prod_{j\in A}\sum_{k_j}|\lambda_{j,{k_j}}|
\frac{|Q_{j,k_j}|^{1+\frac{d+1}{n|A|}}}
{(|x-z_{j,k_j}|+l(Q_{j,k_j}))^{n+\frac{d+1}{|A|}-\frac{\alpha}{|A|}}}
\\
&\times
\prod_{j\in A^c}\sum_{k_j}|\lambda_{j,{k_j}}|
\frac{|Q_{j,k_j}|^{1+\frac{d+1}{n|A|}}}
{(|x-z_{j,k_j}|+l(Q_{j,k_j}))^{n+\frac{d+1}{|A|}}}
\chi_{E_A}\Bigg\|_{L^{q(\cdot)}}.
\end{split}
\end{align}

For $x\in\mathbb R^n$,
denote $\frac{1}{s_j(x)}=\frac{1}{p_j(x)}-\frac{\alpha}{n|A|}$,
$j\in A$. Then
$0<r(x)<\infty$ and
$$
\frac{1}{q(x)}=\sum_{j}^m\frac{1}{p_j(x)}-\frac{\alpha}{n}
=\sum_{j\in A}\frac{1}{s_j(x)}+\sum_{j\in A^c}\frac{1}{p_j(x)}.
$$

For convenience, we denote that
$$
U_A=\sum_{k_j}|\lambda_{j,{k_j}}|
\frac{|Q_{j,k_j}|^{1+\frac{d+1}{n|A|}}}
{(|x-z_{j,k_j}|+l(Q_{j,k_j}))^{n+\frac{d+1}{|A|}-\frac{\alpha}{|A|}}}
$$
and
$$
U_{A^c}=\sum_{k_j}|\lambda_{j,{k_j}}|
\frac{|Q_{j,k_j}|^{1+\frac{d+1}{n|A|}}}
{(|x-z_{j,k_j}|+l(Q_{j,k_j}))^{n+\frac{d+1}{|A|}}}.
$$

Applying (\ref{I22}) and the generalized H\"older inequality
with variable exponents $s_j(\cdot)\;,j\in A$,
$p_j(\cdot)\;,j\in A^c$ and $q(\cdot)$ yield that
\begin{align}\label{I2}
\begin{split}
\|I_2\|_{L^{q(\cdot)}}
&\le
C\sum_{A\subset\{1,\ldots,m\}}\Bigg\|
\prod_{j\in A}U_{A}
\prod_{j\in A^c}U_{A^c}\chi_{E_A}\Bigg\|_{L^{q(\cdot)}}\\
&\le
C\sum_{A\subset\{1,\ldots,m\}}
\left(\prod_{j\in A}\|
U_{A}\chi_{E_A}\|_{L^{s_j(\cdot)}}\right)
\left(\prod_{j\in A^c}\|U_{A^c}
\chi_{E_A}\|_{L^{p_j(\cdot)}}\right).
\end{split}
\end{align}

Denote $\theta=\frac{n+\frac{d+1}{|A|}}{n}$ and we can choose $d$
large enough such that
$\theta p_j^->1$ and $\theta s^-_j>1$.
Notice that $\frac{1}{\theta s_j(x)}=
\frac{1}{\theta p_j(x)}-\frac{\alpha/\theta|A|}{n}$.
By Lemma \ref{s2l3}, we get that
\begin{align}\label{I23}
\begin{split}
\prod_{j\in A}\|U_{A}
\chi_{E_A}\|_{L^{s_j(\cdot)}}
&=\prod_{j=A}\left\|\sum_{k_j}|\lambda_{j,k_j}|
\frac{l(Q_{j,k_j})^{n+\frac{d+1}{|A|}}}
{(|x-z_{j,k_j}|+l(Q_{j,k_j}))^{n+\frac{d+1}
{|A|}-\frac{\alpha}{|A|}}}\right\|_{L^{s_j(\cdot)}}\\
&\le
C
\prod_{j\in A}
\left\|\left(\sum_{k_j}|\lambda_{j,k_j}
|{(M_{\alpha/\theta|A|}{\chi_{Q_{j,k_j}}})^{\theta}}\right)^{1/\theta}
\right\|^\theta_{L^{\theta s_i(\cdot)}}\\
&\le
\prod_{j\in A}\left\|\left(\sum_{k_j}
|\lambda_{j,k_j}|{{\chi_{Q_{j,k_j}
}}}\right)^{\frac{1}{\theta}}
\right\|^\theta_{L^{\theta p_j(\cdot)}}\\
&\leq C\prod_{j\in A}\|f_j\|_{{H}^{p_j(\cdot)}},
\end{split}
\end{align}

where the first inequality follows from the following
claim which can be proved easily:
For any $x\in \mathbb R^n$ and $0\le\alpha<\infty$,
there exists a constant $C$
such that
$$
\frac{r^n}{(r+|x-y|)^{n-\alpha}}\le C(M_\alpha\chi_{Q(y,r)})(x),
$$
where $Q(y,r)$ is a cube centered in $y$ and $r$ its side length.

Repeating the similar argument to (\ref{I23})
with $\alpha=0$, we get that

\begin{align}\label{I24}
\begin{split}
\prod_{j\in A^c}\|U_{A^c}
\chi_{E_A}\|_{L^{p_j(\cdot)}}
&=\prod_{j=A^c}\left\|\sum_{k_j}|\lambda_{j,k_j}|
\frac{l(Q_{j,k_j})^{n+\frac{d+1}{|A|}}}
{(|x-z_{j,k_j}|+l(Q_{j,k_j}))^{n+\frac{d+1}{|A|}}}\right\|_{L^{p_j(\cdot)}}\\
&\le
C
\prod_{j\in A^c}
\left\|\sum_{k_j}|\lambda_{j,k_j}|{(M{\chi_{Q_{j,k_j}}})^{\theta}}
\right\|_{L^{p_i(\cdot)}}\\
&\le
\prod_{j\in A^c}\left\|\left(\sum_{k_j}
|\lambda_{j,k_j}|{{\chi_{Q_{j,k_j}
}}}\right)^{\frac{1}{\theta}}
\right\|^\theta_{L^{\theta p_j(\cdot)}}\\
&\leq C\prod_{j\in A^c}\|f_j\|_{{H}^{p_j(\cdot)}}.
\end{split}
\end{align}
Therefore, for any
$f_j\in{H}^{p_j(\cdot)}\cap{L}^{p_j^-+1}$
by the estimates (\ref{I1}), (\ref{I2}), (\ref{I23}) and (\ref{I24})
then we have
\begin{align}\label{I0}
\|I_\alpha(\vec{f})\|_{L^{q(\cdot)}}\le
C\prod_{j=1}^m\|f_j\|_{{H}^{p_j(\cdot)}}.
\end{align}

From Remark 4.12 in \cite{NS},
we have that
${H}^{p_j(\cdot)}\cap{L}^{\bar p}$
is dense in
${H}^{p_j(\cdot)}$.
Thus,
by the density argument we prove Theorem \ref{s1th1}.

$\hfill\Box$

\section{Commutators of multilinear fractional type operators
on certain Hardy spaces with variable exponents}

In this section, we will study continuity properties of
commutators of multilinear fractional type operators on product of certain
Hardy spaces with variable exponents.
The results are even new for the linear case in the variable exponents
setting.
First we introduce
a new atomic Hardy space with variable exponent
$H^{p(\cdot)}_b$.

\begin{definition}\label{s3de1}
Let b be a locally integrable function and $d\gg 1$. It is said that a
bounded function $a$ is a $(p(\cdot), b, d, \infty)-$atom
if it satisfies\\
\noindent (1) $\mbox{supp}\;a\subset Q = Q(x_0, r)$ for some $r > 0$;\\
\noindent (2) $\|a\|_{\infty}\le \frac{1}{\|\chi_Q\|_{L^{p(\cdot)}}}$;\\
\noindent (3) $\int_{\mathbb R^n} a(x)x^{\alpha}dx
=\int_{\mathbb R^n}a(x)b(x)x^{\alpha}dx=0$\quad for any $|\alpha|\le d$.\\
A temperate distribution f is said to belong to the atomic Hardy space
$H^{p(\cdot)}_b$, if it can
be written as $$f=\sum_{j=1}^\infty\lambda_ja_j,$$ in $\mathcal S'-$sense,
where $a$ is a $(p(\cdot), b, d, \infty)-$atom and
$$
  \mathcal{A}(\{\lambda_j\}_{j=1}^\infty,\{Q_j\}_{j=1}^\infty)
  :=\left\|\left\{\sum_{j}\left(\frac{|\lambda_j|\chi_{Q_j}}
  {\|\chi_{Q_j}\|_{L^{p(\cdot)}}}\right)^{p^-}
\right\}^{\frac{1}{p^-}}\right\|_{L^{p(\cdot)}}<\infty.
$$
Moreover, we define the quasinorm on ${H}_b^{p(\cdot)}$ by
\begin{equation*}
  \|f\|_{{H}_b^{p(\cdot)}}=
  \inf\mathcal{A}(\{\lambda_j\}_{j=1}^\infty,\{Q_j\}_{j=1}^\infty).
\end{equation*}
where the infimum is taken over all admissible expressions
as in $f=\sum_{j=1}^\infty\lambda_ja_j$.
\end{definition}

Obviously, the new atomic Hardy space with variable exponent ${H}_b^{p(\cdot)}$
is a subspace of the Hardy space with variable exponent ${H}^{p(\cdot)}$.
When $p(\cdot)=\mbox{constant}$, the spaces $H^1_b$ and $H^p_b$
first appeared in \cite{P} and \cite{A}. Now we state our second results.

\begin{theorem}\label{s3th1}\quad
Let $0<\alpha<n$
and $b\in BMO$. Suppose that $p_1(\cdot),\ldots,p_m(\cdot)
\in LH\cap \mathcal P^0$ and
$q(\cdot)
\in \mathcal P^0$ be Lebesgue measure functions satisfying
\begin{align}\label{s3c2}
\frac{1}{p_1(x)}+\cdots+\frac{1}{p_m(x)}-\frac{\alpha}{n}=\frac{1}{q(x)},
\quad\quad x\in\mathbb R^n.
\end{align}
Then for any $1\le j\le m$ the operator
$[b,I_\alpha]_j$ can be extended to a bounded operator from
$\prod_{i=1}^m H_{b}^{p_i(\cdot)}$ into $L^{q(\cdot)}$
which satisfies the norm estimate
$$
\|[b,I_\alpha]_j(\vec{f})\|_{L^{q(\cdot)}}
\le C\|b\|_{BMO}\prod_{i=1}^m\|f_i\|_{H_{b}^{p_i(\cdot)}}.
$$
\end{theorem}

%To prove this theorem, we need the two following lemmas.

%\begin{lemma}\cite{TLZ}\label{s3l1}\quad For $b\in L^1_{\rm loc}$, $0<\alpha<nm$ and
%$p_i(\cdot)\in LH\cap\mathcal{P}$, $i=1,2,\cdots,m$. Let $q(\cdot)$ satisfy
%$$\sum_{i=1}^m\frac{1}{p_i(x)}-\frac{\alpha}{n}=\frac{1}{q(x)}<1,\quad x\in\mathbb{R}^n.$$
%Then $[b,I_\alpha]_j$ is bounded from $L^{p_1(\cdot)}\times L^{p_2(\cdot)} \times \cdots\times
%L^{p_m(\cdot)}$ to  $L^{q(\cdot)}$ if and only if $b\in BMO$.
%Furthermore,
%\begin{equation*}
%  \|b\|_{\ast}\sim \|[b,I_\alpha]_j\|_{\mathcal{B}(\prod_{i=1}^mL^{p_i(\cdot)},L^{q(\cdot)})}.
%\end{equation*}
%\end{lemma}

%\begin{lemma}\cite{GK,Tan}\label{s3l2}
%Let $p(\cdot)\in LH\cap \mathcal P^0.$
%Suppose that we are given a sequence of cubes
%$\{Q_j\}_{j=1}^\infty$ and a sequence of non-negative
%$L^1$-functions $\{F_j\}_{j=1}^\infty$.
%Then for $\sum_{j=1}^\infty\chi_{Q_j}F_j \in L^{p(\cdot)}$ we have
%\begin{align*}\label{}
%\left\|\sum_{j=1}^\infty\chi_{Q_j}F_j\right\|_{L^{p(\cdot)}}
%\le C\left\|\sum_{j=1}^\infty\left(\frac{1}{|Q_j|}
%\int_{Q_j}F_j(y)dy\right)\chi_{Q_j}\right\|_{L^{p(\cdot)}}.
%\end{align*}
%\end{lemma}

\noindent\textit{\bf Proof}\quad
The idea of the proof is similar to Theorem \ref{s1th1}.
We only show the differences.
For each
$f_i\in{H}_b^{p_i(\cdot)}$, $i=1,\ldots,m$,
$f_i=\sum_{k_i=1}^\infty\lambda_{i,k_i}a_{i,k_i}$,
where $a_i$ is a $(p_i(\cdot), b, d_i, \infty)-$atom
($d_i\geq d_{p_i(\cdot)}$) and
$$
   \|f_i\|_{{H}_b^{p_i(\cdot)}}=
   \inf\mathcal{A}(\{\lambda_{i,k_i}\}_{i=1}^\infty,
     \{Q_{i,k_i}\}_{i=1}^\infty).
$$

For the decomposition of $f_i$, $i=1,\ldots,m$,
we can write
\begin{align*}\label{}
[b,I_\alpha]_j(\vec{f})(x)=\sum_{k_1}\cdots\sum_{k_m}\lambda_{1,k_1}
\cdots\lambda_{m,k_m}[b,I_\alpha]_j(a_{1,k_1},\ldots,a_{m,k_m})(x)
\end{align*}
in the sense of distributions.

We follow the standard argument in the previous chapter.
Fixed $k_1,\cdots,k_m$, there are two cases for $x\in\mathbb R^n$.\\
Case 1: $x\in Q^\ast_{1,k_1}\cap\cdots\cap Q^\ast_{m,k_m}$.\\
Case 2: $x\in Q^{\ast,c}_{1,k_1}\cup\cdots\cup Q^{\ast,c}_{m,k_m}$.\\

Then we have
\begin{align*}\label{}
&|[b,I_\alpha]_j(\vec{f})(x)|\\
\le&\sum_{k_1}\cdots\sum_{k_m}|\lambda_{1,k_1}|
\cdots|\lambda_{m,k_m}||[b,I_\alpha]_j(a_{1,k_1},\ldots,a_{m,k_m})(x)|
\chi_{Q^\ast_{1,k_1}\cap\cdots\cap Q^\ast_{m,k_m}}(x)\\
&+\sum_{k_1}\cdots\sum_{k_m}|\lambda_{1,k_1}|
\cdots|\lambda_{m,k_m}||[b,I_\alpha]_j(a_{1,k_1},\ldots,a_{m,k_m})(x)|
\chi_{Q^{\ast,c}_{1,k_1}\cup\cdots\cup Q^{\ast,c}_{m,k_m}}(x)\\
=&J_1(x)+J_2(x).
\end{align*}

Now Let us discuss the term $J_1(x)$.
We will show that
\begin{align}\label{J1}
\begin{split}
\|J_1\|_{L^{q(\cdot)}}
\leq C\|b\|_{BMO}\prod_{j=1}^m\|f_j\|_{{H}_b^{p_j(\cdot)}}.
\end{split}
\end{align}
For fixed $k_1,\ldots,k_m$,
assume that $Q^\ast_{1,k_1}\cap\cdots\cap Q^\ast_{m,k_m}\neq 0$, otherwise
there is nothing to prove. Without loss of generality,
suppose that $Q_{1,k_1}$ has the smallest size among all these cubes.
We can pick a cube $G_{k_1,\cdots,k_m}$ such that
$$Q^\ast_{1,k_1}\cap\cdots\cap Q^\ast_{m,k_m}\subset
G_{k_1,\ldots,k_m}\subset G^\ast_{k_1,\ldots,k_m}\subset
Q^{\ast\ast}_{1,k_1}\cap\cdots\cap Q^{\ast\ast}_{m,k_m}$$
and $|G_{k_1,\ldots,k_m}|\ge C|Q_{1,k_1}|$.

Denote $M(x):=|[b,I_\alpha]_j(a_{1,k_1},\ldots,a_{m,k_m})(x)|
\chi_{Q_{1,k_1}^\ast\cap\cdots\cap Q_{m,k_m}^\ast}(x)$.
Obviously, $$\mbox{supp} M(x)\subset
{Q_{1,k_1}^\ast\cap\cdots\cap Q_{m,k_m}^\ast}\subset
G_{k_1,\cdots,k_m}.$$

By \cite[Theorem 1.1]{TLZ},
we have that $[b,I_\alpha]_j$ maps $L^{r_1}\times L^\infty\times\cdots
\times L^\infty$ into $L^{s}$ ($s>1$ and
$\frac{1}{s}=\frac{1}{r_1}-\frac{\alpha}{n}$).

Then we
get that
\begin{align*}\label{}
\begin{split}
\|M\|_{L^{s}}&\le
\|[b,I_\alpha]_j(a_{1,k_1},\ldots,a_{m,k_m})\|_{L^{s}}\\
&\le
C\|b\|_{BMO}\|a_{1,k_1}\|_{L^{r_1}}\|a_{2,k_2}\|_{L^\infty}
\cdots\|a_{m,k_m}\|_{L^\infty}\\
&\le C\|b\|_{BMO}\prod_{i=1}^m\|\chi_{Q_{i,{k_i}}}\|^{-1}_{L^{p_i(\cdot)}}
|G_{k_1,\ldots,k_m}|^{\frac{1}{r_1}}.
\end{split}
\end{align*}

Repeating similar argument to the proof of
(\ref{I1}) in Theorem \ref{s1th1}, we can
obtain the desired result (\ref{J1}).

Secondly, we consider the estimate of $J_2$.
When $x\in\chi_{Q^{\ast,c}_{1,k_1}\cup\cdots\cup Q^{\ast,c}_{m,k_m}}(x)$,
the $L^{q(\cdot)}$ norm of $J_2$ is controlled by
\begin{align*}\label{}
&\left\|\sum_{k_1}\cdots\sum_{k_m}|\lambda_{1,k_1}|
\cdots|\lambda_{m,k_m}|
(b(x)-b_{Q_{j,k_j}})I_\alpha
(a_{1,k_1},\ldots,a_{m,k_m})\right\|_{L^{q(\cdot)}}\\
&+\left\|\sum_{k_1}\cdots\sum_{k_m}|\lambda_{1,k_1}|
\cdots|\lambda_{m,k_m}|
|I_\alpha(a_{1,k_1},\ldots,
(b-b_{Q_{j,k_j}})a_{j,k_j},\ldots,a_{m,k_m})\right\|_{L^{q(\cdot)}}\\
&=:\|J_{21}\|_{L^{q(\cdot)}}+\|J_{22}\|_{L^{q(\cdot)}}.
\end{align*}

Let A be a nonempty
subset of $\{1,\ldots,m\}$, and we denote the cardinality of
$A$ by $|A|$, then $1\le|A|\le m$. Let $A^c=\{1,\ldots,m\}\backslash
A.$ If $A=\{1,\ldots,m\}$, we define
$$(\cap_{i\in A}Q_{i,k_i}^{\ast,c})
\cap (\cap_{i\in A^c}Q_{i,k_i}^{\ast,c})=
\cap_{i\in A}Q_{j,k_i}^{\ast,c},$$
then
$$Q_{1,k_1}^{\ast,c}\cup\cdots\cup Q_{m,k_m}^{\ast,c}
=\cup_{A\subset\{1,\ldots,m\}}((\cap_{i\in A}Q_{i,k_i}^{\ast,c})
\cap (\cap_{i\in A^c}Q_{i,k_i}^{\ast})).$$

Set $E_A=(\cap_{i\in A}Q_{i,k_i}^{\ast,c})
\cap (\cap_{i\in A^c}Q_{i,k_i}^{\ast}).$
For fixed $A$, assume that $Q_{\tilde i,k_{\tilde i}}$
is the smallest cubes in the set of cubes $Q_{i,k_i},
i\in A$.
Let $z_{\tilde i,k_{\tilde i}}$ is the center of the cube
$Q_{\tilde i,k_{\tilde i}}$.

From Definition \ref{s3de1} and (\ref{I20}), we can easily get that
\begin{align}\label{J20}
\begin{split}
&|I_\alpha(a_{1,k_1},\ldots,a_{m,k_m})(x)|\\
\le&
C\prod_{i\in A}\frac{|Q_{i,k_i}|^{1+\frac{d+1}{n|A|}}}
{\|\chi_{Q_{i,{k_i}}}\|_{L^{p_i(\cdot)}}
(|x-z_{i,k_i}|+l(Q_{i,k_i}))^{n+\frac{d+1}{|A|}-\frac{\alpha}{|A|}}}\\
&\times \prod_{i\in A^c}\frac{|Q_{i,k_i}|^{1+\frac{d+1}{n|A|}}}
{\|\chi_{Q_{i,{k_i}}}\|_{L^{p_i(\cdot)}}
(|x-z_{i,k_i}|+l(Q_{i,k_i}))^{n+\frac{d+1}{|A|}}}\\
=:&C\prod_{i\in A}U_i \prod_{i\in A^c}V_i
\end{split}
\end{align}
for all $x\in E_A$.

For $x\in\mathbb R^n$,
denote $\frac{1}{s_i(x)}=\frac{1}{p_i(x)}-\frac{\alpha}{n|A|}$,
$i\in A$. Then
$0<r(x)<\infty$ and
$$
\frac{1}{q(x)}=\sum_{i}^m\frac{1}{p_i(x)}-\frac{\alpha}{n}
=\sum_{i\in A}\frac{1}{s_i(x)}+\sum_{i\in A^c}\frac{1}{p_i(x)}.
$$

We will discuss in two case to estimate $J_2$.
When $j\in A$, by generalized H\"older's inequality
in variable Lebesgue spaces we have
\begin{align*}\label{}
\|J_{21}\|_{L^{q(\cdot)}}
\le&\left\|\sum_{k_1}\cdots\sum_{k_m}|\lambda_{1,k_1}|
\cdots|\lambda_{m,k_m}|
(b(x)-b_{Q_{j,k_j}})
\prod_{i\in A}U_i\prod_{i\in A^c}V_i\chi_{E_A}
\right\|_{L^{q(\cdot)}}\\
\le&
\prod_{i\in A;i\neq j}\left\|\sum_{k_i}|\lambda_{i,k_i}|
U_i\chi_{E_A}\right\|_{L^{s_i(\cdot)}}
\left\|\sum_{k_j}|\lambda_{j,k_j}|(b(x)-b_{Q_{j,k_j}})
U_j\chi_{E_A}\right\|_{L^{s_j(\cdot)}}\\
&\times\prod_{i\in A^c}\left\|\sum_{k_i}|
\lambda_{i,k_i}|V_i\chi_{E_A}\right\|_{L^{p_i(\cdot)}}.
\end{align*}

We follow the similar argument in the previous chapter agian.
Denote $\theta=\frac{n+\frac{d+1}{|A|}}{n}$ and we can choose $d$
large enough such that
$\theta p_j^->1$ and $\theta s^-_j>1$.
Notice that $\frac{1}{\theta s_j(x)}=
\frac{1}{\theta p_j(x)}-\frac{\alpha/\theta|A|}{n}$.
By Lemma \ref{s2l3}, we get that
\begin{align}\label{J23}
\begin{split}
&\prod_{i\in A;i\neq j}\left\|\sum_{k_i}|\lambda_{i,k_i}|
U_i\chi_{E_A}\right\|_{L^{s_i(\cdot)}}
\le
C
\prod_{i\in A;i\neq j}
\left\|\left(\sum_{k_i}|\lambda_{i,k_i}
|{(M_{\alpha/\theta|A|}{\chi_{Q_{i,k_i}}})^{\theta}}\right)^{1/\theta}
\right\|^\theta_{L^{\theta s_i(\cdot)}}\\
&\le
\prod_{i\in A;i\neq j}\left\|\left(\sum_{k_i}
|\lambda_{i,k_i}|{{\chi_{Q_{i,k_i}
}}}\right)^{\frac{1}{\theta}}
\right\|^\theta_{L^{\theta p_i(\cdot)}}
\leq C\prod_{i\in A;i\neq j}\|f_j\|_{{H}^{p_j(\cdot)}}.
\end{split}
\end{align}

Repeating the similar but easier argument to (\ref{I23})
with $\alpha=0$, we get that

\begin{align}\label{J24}
\begin{split}
\prod_{i\in A^c}\left\|\sum_{k_i}|
\lambda_{i,k_i}|V_i\chi_{E_A}\right\|_{L^{p_i(\cdot)}}
\leq C\prod_{i\in A^c}\|f_i\|_{{H}^{p_i(\cdot)}}.
\end{split}
\end{align}

Next we need to prove
\begin{align}\label{J25}
\begin{split}
\left\|\sum_{k_j}|\lambda_{j,k_j}|(b(x)-b_{Q_{j,k_j}})
U_j\chi_{E_A}\right\|_{L^{s_j(\cdot)}}
\le C\|f_j\|_{{H}^{p_j(\cdot)}}.
\end{split}
\end{align}

For any constants $s,s_i>1,\;i=1,2,3,$
denote that $\frac{1}{s}=\frac{1}{s_1}+\frac{1}{s_2}
=\frac{1}{s_1}+\frac{1}{r_2}-\frac{\alpha}{n|A|}=:\frac{1}{s_3}-\frac
{\alpha}{n|A|}$.

Applying H\"older's inequality yields that
\begin{align*}\label{}
\begin{split}
\left\|(b(x)-b_{Q_{j,k_j}})
U_j\chi_{E_A}\right\|_{L^s}
&\le C\left\|(b(x)-b_{Q_{j,k_j}})
U_j\chi_{E_A}\right\|_{L^{s_1}}
\left\||
U_j\chi_{E_A}\right\|_{L^{s_2}}\\
&\le C\|b\|_{BMO}|Q|^{1/{s_1}}
\left\||
(M_{\alpha/\theta|A|}{\chi_{Q_{i,k_i}}})^{\theta}
\right\|_{L^{s_2}}\\
&\le C\|b\|_{BMO}|Q|^{1/{s_1}}
|Q|^{1/{r_2}}\\
&\le C\|b\|_{BMO}|Q|^{1/{s_3}},
\end{split}
\end{align*}
where the second estimate comes from
the John-Nirenberg inequality.

Now we resort to the proof of (\ref{I1})
to have the desired results (\ref{J25}).

Therefore, when $j\in A$,
by the estimates (\ref{J23}), (\ref{J24}) and (\ref{J25})
we have
\begin{align}\label{J21}
\|J_{21}\|_{L^{q(\cdot)}}\le
C\prod_{i=1}^m\|f_i\|_{{H}_b^{p_j(\cdot)}}.
\end{align}

When $j\in A^c$, we similarly have
\begin{align*}\label{}
\|J_{21}\|_{L^{q(\cdot)}}
\le&
\prod_{i\in A}\left\|\sum_{k_i}|\lambda_{i,k_i}|
U_i\chi_{E_A}\right\|_{L^{s_i(\cdot)}}
\left\|\sum_{k_j}|\lambda_{j,k_j}|(b(x)-b_{Q_{j,k_j}})
V_j\chi_{E_A}\right\|_{L^{p_j(\cdot)}}\\
&\times\prod_{i\in A^c;i\neq j}\left\|\sum_{k_i}|
\lambda_{i,k_i}|V_i\chi_{E_A}\right\|_{L^{p_i(\cdot)}}.
\end{align*}

We only need to observe that
\begin{align*}\label{}
\begin{split}
\left\|(b(x)-b_{Q_{j,k_j}})
V_j\chi_{E_A}\right\|_{L^s}
&\le C\left\|(b(x)-b_{Q_{j,k_j}})
U_j\chi_{E_A}\right\|_{L^{s}}
\left\||
U_j\chi_{E_A}\right\|_{L^{\infty}}\\
&\le C\|b\|_{BMO}|Q|^{1/{s}}
\left\||
(M{\chi_{Q_{i,k_i}}})^{\theta}
\right\|_{L^{\infty}}\\
&\le C\|b\|_{BMO}|Q|^{1/{s}}.
\end{split}
\end{align*}

Similarly, we can get (\ref{J21}).

Let us estimate $J_{22}$.
Observe that
$\|a_{i,k_i}\|_{\infty}\le \frac{1}{\|\chi_{Q_{i,k_i}}\|_{L^{p_i(\cdot)}}}$ and
$\int_{\mathbb R^n} a_{\tilde i,k_{\tilde i}}(x)x^{\alpha}dx
=\int_{\mathbb R^n}a_{\tilde i,k_{\tilde i}}(x)b(x)x^{\alpha}dx=0$\quad for any $|\alpha|\le d$.

As the argument for (\ref{I20}),
we similarly have that
\begin{align*}\label{}
\begin{split}
&|I_\alpha(a_{1,k_1},\ldots,
(b-b_{Q_{j,k_j}})a_{j,k_j},\ldots,a_{m,k_m})|\\
\le&
C\|b\|_{BMO}\prod_{i\in A}\frac{|Q_{i,k_i}|^{1+\frac{d+1}{n|A|}}}
{\|\chi_{Q_{i,{k_i}}}\|_{L^{p_i(\cdot)}}
(|x-z_{i,k_i}|+l(Q_{i,k_i}))^{n+\frac{d+1}{|A|}-\frac{\alpha}{|A|}}}\\
&\times \prod_{i\in A^c}\frac{|Q_{i,k_i}|^{1+\frac{d+1}{n|A|}}}
{\|\chi_{Q_{i,{k_i}}}\|_{L^{p_i(\cdot)}}
(|x-z_{i,k_i}|+l(Q_{i,k_i}))^{n+\frac{d+1}{|A|}}}
\end{split}
\end{align*}
for all $x\in E_A$.
The rest of the proof is same as the above.
Then we have
\begin{align}\label{J22}
\|J_{22}\|_{L^{q(\cdot)}}\le
C\prod_{i=1}^m\|f_i\|_{{H}_b^{p_j(\cdot)}}.
\end{align}

In conclusion, combing \ref{J1}, \ref{J21} and \ref{J22},
we obtain that
\begin{align*}\label{}
\|[b,I_\alpha]_j(\vec{f})\|_{L^{q(\cdot)}}
\le C\|b\|_{BMO}\prod_{i=1}^m\|f_i\|_{H_{b}^{p_i(\cdot)}}.
\end{align*}

We completed
the proof of Theorem \ref{s3th1}.
$\hfill\Box$

{\bf Acknowledgments.}
The author wish to express his heartfelt thanks to the anonymous referees for careful reading.
The project is sponsored by Natural Science Foundation of Jiangsu Province of China (grant no. BK20180734), Natural Science Research of Jiangsu Higher Education Institutions of China (grant no. 18KJB110022) and Nanjing University of Posts and Telecommunications Science Foundation (grant no. NY219114).

\end{document}